\documentclass[times,doublespace]{oupau}
\usepackage{amsmath,amssymb,graphicx,amscd,hyperref,amsfonts,chemarr,graphicx,psfrag,fancybox,qsymbols,indentfirst}
\usepackage{mathrsfs}
 \usepackage[height=190mm,width=130mm]{geometry}
\usepackage{dsfont}
\usepackage[latin1]{inputenc}

\def \sous#1#2{\mathrel{\mathop{\kern 0pt#1}\limits_{#2}}}
\def \sur#1#2{\mathrel{\mathop{\kern 0pt#1}\limits^{#2}}}

\def \el{\sur{=}{(d)}}

\def \be{\begin{eqnarray*}}
\def \ee{\end{eqnarray*}}
\def \ben{\begin{eqnarray}}
\def \een{\end{eqnarray}}

\def\wt{\widetilde}

\def\sn{^{(N)}}

\def\AArm{\fam0 }
\def\AAk#1#2{\setbox\AAbo=\hbox{#2}\AAdi=\wd\AAbo\kern#1\AAdi{}}%
\def\AAr#1#2#3{\setbox\AAbo=\hbox{#2}\AAdi=\ht\AAbo\raise#1\AAdi\hbox{#3}
}%

\def \sur#1#2{\mathrel{\mathop{\kern 0pt#1}\limits^{#2}}}
\def\BBr{{\AArm I\!R}}%

\def \a{{\tt a}}
\def \b{{\tt b}}

\def \w{{\tt w}}
\def \u{{\tt u}}

\def\videbox{\mathbin{\vbox{\hrule\hbox{\vrule height1ex \kern.5em\vrule height1ex}\hrule}}}

\def \build#1#2#3{\mathrel{\mathop{\kern 0pt#1}\limits_{#2}^{#3}}}

\def\proof{\noindent{\bf Proof}\hskip10pt}
\def\QED{\hfill\vrule height 1.5ex width 1.4ex depth -.1ex \vskip20pt}

\begin{document}
\runningheads{F. Gamboa and A. Rouault}{Large deviations and sum rules}
\title{Large deviations for random spectral measures and sum rules}

\author{Fabrice Gamboa\affil{a},Alain Rouault\affil{b}}

\address{\affilnum{a}Institut de Math\'ematiques de Toulouse.  
Universit\'e Paul Sabatier,                                                               
F-31062 Toulouse Cedex 9
\affilnum{b}LMV  B\^atiment Fermat,
 Universit\'e
Versailles-Saint-Quentin,
F-78035 Versailles Cedex}
\corraddr{LMV  B\^atiment Fermat,
 Universit\'e
Versailles-Saint-Quentin,
F-78035 Versailles Cedex. E-mail: Alain.Rouault@math.uvsq.fr}
\begin{abstract}
We prove a Large Deviation Principle for the random spectral measure associated to the pair $(H_N, e)$ where 
$H_N$ is sampled in the GUE and $e$ is a fixed unit vector (and more generally in the $\beta$ extension of this model). The rate function consists of two parts. The contribution of the absolutely continuous part of the measure is the reversed Kullback information with respect to the semicircle distribution and the contribution of the singular part is connected to the rate function of the extreme eigenvalue in the GUE. This method is also applied to the Laguerre and Jacobi ensembles, but in those cases the expression of the rate function is not explicit. 
\end{abstract}
\received{}

\newtheorem{remark}{Remark}
\newtheorem{conj}{Conjecture}

\def\proof{\noindent{\bf Proof}\hskip10pt}
\def\QED{\hfill\vrule height 1.5ex width 1.4ex depth -.1ex \vskip20pt}
%

\maketitle
\section{Introduction}
\label{Sintno}
The aim of this paper is to study the asymptotic behavior of spectral measures in some classical self adjoint random matrix models. To begin with, let us first clarify what we mean with spectral measure of a pair in the case of unitary operators and  recall some asymptotic results in this case.
\\
\\
Let  $U$ be a unitary operator in a Hilbert space $\mathcal H$, and 
$e$ be a unit cyclic vector (the span generated by the iterates
 $(U^nx)$ is $\mathcal H$).  The spectral measure associated with  the pair $(U,e)$ plays an important role and will be one of the object studied here.  This measure is the unique probability measure (p.m.)
 $\mu$ on the unit circle $\mathbb T$ such that
\[\langle e, U^n e\rangle = \int_{\mathbb T} z^n d\mu(z) \ \ \ (n\geq 1)\,.\]
This measure is a unitary invariant for the pair $(U,e)$.
Assume further that dim ${\mathcal H} =N$ and that $e_1$ the  first vector of the canonical basis is cyclic for $U$. Let $\lambda_1, \dots, \lambda_N$ be the eigenvalues of $U$ (all lying on $\mathbb T$),  
and let $\psi_1, \dots, \psi_N$ be a system of unit eigenvectors. 
The spectral measure is then
\ben
\label{masterdef}
\mu_\w\sn = \sum_{k=1}^N \pi_k\delta_{\lambda_k}\,\een
with $\pi_k := |\langle \psi_k, e_1\rangle|^2 \ \ k=1, \dots, N$. Notice that given $\lambda_k$, the vector $\psi_k$ is determined up to a phase, but the number $\pi_k$ is completely determined. 
To avoid confusion, 
we put an index $\w$ (for weight) to distinguish this measure from the classical empirical spectral distribution (ESD) defined by
\ben\mu_\u\sn = \frac{1}{N} \sum_{k=1}^N \delta_{\lambda_k}\,.\een
When $U$ is uniformly sampled from  $\mathbb U(N)$ (the unitary group of order $N$)  with the Haar distribution, it is well known that  the joint distribution of $(\lambda_1, \dots, \lambda_N)$ has a density proportional to
\[|\Delta(\lambda_1, \dots, \lambda_N)|^2\]
where $\Delta$ is the Vandermonde determinant (see for example 
\cite{MehtaB}). Furthermore, $e_1$ is almost surely (a.s.) cyclic and  $(\pi_1, \dots, \pi_N)$ is independent of $(\lambda_1,\dots, \lambda_N)$.  Moreover,  $(\pi_1, \dots, \pi_N)$ is uniformly distributed on the simplex ${\mathcal S}_n = \{(\pi_1, \dots, \pi_N) : \pi_k > 0, (k=1, \dots , N), \!\ \pi_1 + \dots \pi_N = 1 \}$. As $N$ tends to infinity, both sequences of random measures $(\mu_\w\sn)$  and $(\mu_\u\sn)$   converge weakly to the equilibrium measure, i.e. the uniform distribution on $\mathbb T$.
In a previous work (\cite{FABALOZ}), 
we proved  that the sequence  $(\mu_\w\sn)$ satisfies a Large Deviation Principle  (denoted hereafter LDP), with speed $N$ and good rate function given by reversed Kullback entropy  with respect to the equilibrium measure.
Notice that there is a quite important difference in the large deviation behaviour of  $(\mu_\w\sn)$ and $(\mu_\u\sn)$.  Indeed, this last sequence of probability measures (p.ms.) satisfies a LDP with speed $N^2$ and with a rate function connected to the Voiculescu entropy 
(see for example \cite{HiaiIHP}). 
To show a LDP for $(\mu_\w\sn)$ one may think of two kinds of proof. The first one, which could be called the direct way, uses the representation (\ref{masterdef}) \cite{FABALOZ}. Besides, it is possible to code a measure $\mu$  on $\mathbb T$ by the system of its Verblunsky (or Schur) coefficients, via the Favard theorem \cite{Simon1};  they are also the canonical moments of $\mu$ (see \cite{DS} for the definition). The second method uses this coding.  It turns out that, under the Haar distribution,  the canonical moments $(c_1\sn, \dots c_N\sn)$ of $\mu_\w\sn$ are independent random variables (r.vs.) with explicit distribution depending on $N$. It is then possible in a first step to check the LDP on these variables and in a second step to lift the LDP and the rate function on the space of measures \cite{LozEJP}.

The precise form of the rate function can be explained, in the first method by the Dirichlet weighting of the random measure, and in the second method by the Szeg\"o formula, which enters in the class of the so-called sum rules. The same thing can be done for the Jacobi ensemble with the arcsine distribution (on $[0,1]$ or on $[-2, 2]$) playing the role of the uniform distribution on $\mathbb T$. 

In this paper we will focus on  models of self-adjoint matrices and their extensions. If $H$ is a self-adjoint bounded operator in a Hilbert space ${\mathcal H}$ and $e$ a cyclic vector, the spectral measure is the unique p.m. $\mu$ on $\mathbb R$ such that
\[\langle e, H^n e\rangle = \int_{\mathbb R} x^n d\mu(x) \ \ \ (n\geq 1)\,.\]
Here also, $\mu$ is an unitary invariant for the pair $(H,e)$. Another invariant is the tridiagonal reduction recalled in Section \ref{Sjaja}.
If dim ${\mathcal H} =N$ and $e_1$ is cyclic for $H$, the spectral measure is
\ben
\label{remuw}
\mu_\w\sn  = \sum_{k=1}^N \pi_k\delta_{\lambda_k}\,\een
with the same notation as above, except that, now, the eigenvalues are real.

The first two models leads to an eigenvalue distribution that is not
almost surely (a.s.) supported by a fixed compact set.
We will first study the $\beta$-Hermite ensemble. It is a family extending the Gaussian ensembles (GOE, GUE, GSE).  The second model considered is the $\beta$-Laguerre ensemble that generalizes Wishart matrices. In both cases,
we could expect that the sequence  $(\mu_\w\sn)$ satisfies a LDP with speed $N$ and with a rate function given by the reversed Kullback entropy with respect to the limit distribution 
(respectively semicircle and Marchenko-Pastur distributions). Actually the difference with the unitary case comes from the problem of support. 
We  prove results of the same flavour that those we previously obtained in the unitary case, 
but with an extra contribution in the rate function due to the singular part of measures.  The third model studied is the $\beta$-Jacobi ensemble in which the eigenvalues are confined in a compact set.

The paper is organized as follows.  
The next section is devoted to the introduction of  notation and models : topology on space of moments and 
 real matrix models that we will study later. 
In Section \ref{Sjaja}, we discuss  some relationships between the random spectral measures and coefficients appearing in the construction of the associated random orthogonal polynomials. The LDP for real matrix models are studied in last two sections. The case of the $\beta$-Hermite ensemble is completely tackled in Section \ref{sHerm}. Surprisingly, we manage  to compute  explicitely  the rate function, with the help of a convenient sum rule. The $\beta$-Laguerre and $\beta$-Jacobi ensembles are studied in Section \ref{sLag} and \ref{sJac}. Here, the rate functions are not so explicit. 
All useful distributions we work with are defined in Section \ref{appendixb}. After posting a previous version of this paper on arxiv we have been aware of a paper of Dette and Nagel (see \cite{De118}) stating CLTs for moments of the random spectral measures studied here. 
\section{Notation and models}
\subsection{Topology on moments spaces}
Let ${\mathcal M}^1$ be the set of all p.ms. on $\BBr$ and let  
 ${\mathcal M}_m^1$ be the subset consisting of p.ms. on $\BBr$ having all their  moments finite.
  For $\mu \in  {\mathcal M}_m^1$ we set 
\[m_k(\mu) = \int_{\mathbb R} x^k d\mu(x) \ , \ \ k \geq 1\,,\]
and ${\mathbf m}(\mu)= \big(m_k (\mu)\big)_{k \geq 1}$.
As it is classical in moment problems, we consider the set ${\mathcal M}_m^1$ 
 as a subset of $\mathbb R[[X]]$, (the set of formal series with real coefficients), or equivalently as a subset of the set of linear forms on the space $\mathbb R [X]$ of polynomials with real coefficients, or eventually as a subset of ${\mathbb R}^{\mathbb N}$. We may identify $\mu$  either with
\begin{itemize}
\item The formal series $\sum_{n=0}^\infty m_n(\mu) X^n$,
\item The 
 linear form on $\mathbb R[X] : P(X) \mapsto \int_{\mathbb R} P(x) d\mu(x)$
\item The sequence ${\mathbf m}(\mu)$.
\end{itemize}
We endow  ${\mathcal M}_m^1$ with the distance of convergence of moments:
\ben
\label{defdistmunu}
d(\mu, \nu) = \sum_{k=1}^\infty 2^{-k} \frac{|m_k(\mu) - m_k (\nu)|}{1 +|m_k(\mu) - m_k (\nu)|}\,.
\een
If ${\mathcal M}_{m,d}^1$ denotes the subset of ${\mathcal M}_m^1$ consisting in all p.ms. determined by their moments, then the mapping $\mathbf m$ is injective and continuous from  ${\mathcal M}_{m,d}^1$ to ${\mathbb R}^{\mathbb N}$. Notice that the topology on ${\mathcal M}_m^1$ used here is quite different from the usual weak topology. Indeed, on the one hand convergence on bounded continuous function does not imply convergence of the moments. On the other hand it is not possible to approximate uniformly on $\mathbb{R}$ a bounded continuous function by a polynomial.\\
\\
In the next subsections we recall some classical ensembles of random matrices. We refer to \cite{agz} for a complete overview on this topic.
\subsection{$\beta$-Hermite ensemble}
Let us begin with Gaussian matrix models and their extensions. 
\begin{itemize}
\item GOE(N) The diagonal entries are independent and ${\mathcal N}(0; 2/N)$ distributed and the non diagonal entries are independent up to symmetry and ${\mathcal N}(0; 1/N)$ distributed.
The joint density on $\BBr^N$ of the eigenvalues is proportional to
\[\Delta(\lambda_1, \dots, \lambda_N) \exp - \frac{N}{4} \sum_j \lambda_j^2\,.\]  
The matrix of eigenvectors is orthogonal, so its first line is uniformly distributed on the $N$-dimensional sphere, i.e. the vector $(\pi_1, \dots, \pi_N)$ has the distribution Dir$_N(1/2)$.
\item GUE($N$) The diagonal entries are independent and ${\mathcal N}(0; 1/N)$ distributed and the non diagonal entries are independent up to symmetry and distributed as ${\mathcal N}(0; 1/2N) + \sqrt{-1}\!\ {\mathcal N}(0; 1/2N)$ where both normal variables are independent. 
The joint density of the eigenvalues is proportional to
\[\Delta(\lambda_1, \dots, \lambda_N)^2 \exp - \frac{N}{2} \sum_j \lambda_j^2.\] 
The matrix of eigenvectors is unitary, so the first line is uniformly distributed on the $N$-dimensional (complex) sphere, i.e. the vector $(\pi_1, \dots, \pi_N)$ 
has the distribution Dir$_N (1)$. 

If $M$ is sampled from the GOE($N$) or GUE($N$),  $e_1$ is a.s. cyclic, the eigenvalues are a.s. distinct and then 
we will consider the (random) spectral measure $\mu_\w\sn$ given by (\ref{remuw}).

We do not 
recall the definition of
the symplectic ensemble GSE($N$). Nevertheless, some of the previous objects may also be defined in this context.

\item More generally, it is now classical to consider a parameter $\beta = 2 \beta' >0$, and a density in $\mathbb R^N$ proportional to
\[|\Delta(\lambda_1, \dots, \lambda_N)|^{\beta} \exp -\frac{N\beta}{4} \sum_j \lambda_j^2\,.\]
This expression extends the above formulas so that $\beta= 1$ for the GOE, $\beta=2$ for the GUE and $\beta=4$ for the GSE.
It is often called a Coulomb gas model and $(\lambda_1, \dots, \lambda_N)$ are called charges.

Dumitriu and Edelman (\cite{DuEdel2} Theorem 2.12) found a matrix model for this distribution, i.e. a random real symmetric matrix whose eigenvalues follows the above distribution. Moreover they proved that the corresponding vector $(\pi_1, \dots, \pi_N)$ is independent of the eigenvalues and  Dir$_N(\beta')$ distributed. A specific description of the matrix will be given in the next section.   
 
When $N \rightarrow \infty$, it is known that $(\mu_\u\sn)$  converges weakly to the semicircle distribution, and satisfies a LDP with speed $N^2$ and with a rate function connected to the Voiculescu entropy.
\end{itemize} 
\subsection{$\beta$-Laguerre and $\beta$-Jacobi ensembles}
\begin{itemize}
\item The classical Wishart real ensemble is formed by   
$W= G\ ^tG$ with $G$ a $m\times N$ matrix with independent ${\mathcal N}(0, 2/N)$ entries.
The joint density of eigenvalues is proportional to
\[ |\Delta(\lambda)| \prod_{j=1}^m \lambda_j^{\frac{1}{2}(N-m+1)-1} \exp -\frac{N}{4}\sum_{j=1}^m \lambda_j\]
and the distribution of weights $(\pi_1, \dots, \pi_m)$  is Dir$_N(1/2)$.

 This eigenvalues distribution is classicaly extended to the $\beta$-Laguerre distribution of charges, with density proportional to :
\[ |\Delta(\lambda)|^\beta \prod_{j=1}^m \lambda_j^{\beta'(N-m+1)-1} \exp -\frac{N\beta}{4}\sum_{j=1}^m \lambda_j\,.\]
For this case, Dumitriu and Edelman (\cite{DuEdel2} Theorem 3.4) also gave a (real symmetric) matrix model and proved that the vector of 
  weights $(\pi_1, \dots, \pi_m)$ is also independent of the eigenvalues and is Dir$_N(\beta')$ distributed.
\item The $J\beta E(N; \a, \b)$ ensemble (with $\a > -1, \b> -1$) has been defined to extend the MANOVA ensemble known in statistics for $\beta =1$ and $\beta=2$. It is defined by  a density of $N$ charges on $[-2, 2]$
\[|\Delta(x_1, \dots, x_N)|^\beta \prod_{j=1}^N (2 -x_j)^\a (2+ x_j)^\b\,.\] 
Killip and Nenciu (\cite{Killip1}) found a matrix (real symmetric) model and proved that the corresponding vector of weights is again independent of the eigenvalues and  Dir$_N(\beta')$ distributed.
A variant is the $\widehat{J\beta E}(N, \a, \b)$ ensemble where the charges  are distributed on $[0, 1]$ according to a density proportional to
\[|\Delta(x_1, \dots, x_N)|^\beta \prod_{j=1}^N x_j^\a (1- x_j)^\b\,,\] 
In the matrix model,  the weights have the same properties as above.
\end{itemize}
\section{Tridiagonal representations}
\label{Sjaja}
\subsection{Spectral map}
In this section, we will describe the Jacobi mapping between tridiagonal matrices and spectral measures.  This mapping will be one of the key tools for our large deviations results. We consider finite size matrices  corresponding to  measures supported by a finite number of points and semi-infinite matrices corresponding to measures with bounded infinite support. The material of this section is largely borrowed from \cite{Simon3},
\cite{simon2006cmf} \cite{Simonself}. 
 
If $\mu$ is a probability measure with a finite support consisting of $N$ points 
the orthonormal polynomials (with positive leading coefficients) obtained by Gram-Schmidt procedure from the sequence $1, x, x^2, \dots, x^{N-1}$ satisfy the recurrence relation
\[xp_n (x) = a_{n}p_{n+1}(x) + b_{n}p_n(x) + a_{n-1} p_{n-1}(x)\]
for $n \leq N-1$, where $a_n > 0$ for those $n$. In the basis $\{p_0, p_1, \dots, p_{n-1}\}$, the linear transformation $x \mapsto xf(x)$ in $L^2(d\mu)$ is represented by the matrix
\ben
\label{favardfini}
J_\mu = \begin{pmatrix} b_0&a_0 &0&\dots&0\\
a_0&b_1  &a_1&\ddots&\vdots\\
0&\ddots &\ddots&\ddots&0\\
\vdots&\ddots&a_{N-3}&b_{N-2}&a_{N-2}\\
0&\dots&0&a_{N-2}&b_{N-1}
\end{pmatrix}
\een
So, measures supported by $N$ points lead to Jacobi matrices, i.e. $N\times N$ symmetric tridiagonal matrices with subdiagonal positive terms ; in fact, there
is a one-to-one correspondence between them. Given  such a Jacobi matrix $J$, $e_1$ is cyclic and if $\mu$ is 
the spectral measure associated to the pair $(J, e_1)$, then 
 $J$
represents the multiplication by $x$ in the basis of orthonormal polynomials associated to $\mu$ and $J=J_\mu$.

More generally, if $\mu$ is a p.m. on $\BBr$, with bounded infinite support, we may apply the same Gram-Schmidt process and consider the associated semi-infinite Jacobi matrix: 
\ben
\label{favardinfini}
J_\mu = \begin{pmatrix} b_0&a_0 &0&0&\dots\\
a_0&b_1  &a_1&0&\dots\\
0&a_1 &b_2&a_3&\dots\\
\dots&\dots&\dots&\dots&\dots
\end{pmatrix}
\een
Notice that  again we have $a_k > 0$ for every $k$. The mapping $\mu \mapsto J_\mu$ (which we call Jacobi mapping) 
 is a one to one correspondence between p.ms on $\mathbb R$ having compact infinite support and this kind of tridiagonal matrices 
with  $\sup_n(|a_n| +
|b_n|) < \infty$. This result is 
sometimes called Favard's theorem (see \cite{Simon3} p.432).

Furthermore, a compactly supported p.m. $\mu$  is  completely determined through the knowledge of all  its moments $m_k(\mu)$  for $k \geq 1$. So, an inversion formula for the Jacobi mapping may be performed by using $J_\mu$ to compute the moments of $\mu$ (see for example \cite{Simonself}). 
Actually, it is a recursive procedure and it is possible to connect successive moments with sucessive sections of the matrix. For a general Jacobi semi-infinite (resp. $N\times N$) matrix $A$, let $A^{[j]}$ for $j \geq 1$ (resp. for $j \leq N$) the left top submatrix of $A$. 
It is known from \cite{Simonself} formula (5.37), that if $A$ is semi-infinite, we have the identity
\ben
\langle e_1, A^k e_1\rangle &=& \langle e_1, \big(A^{[j]}\big)^{k}e_1\rangle \ , \ k=1, \dots, 2j-1.
\een
It is straightforward that this formula holds true when $A$ is a Jacobi $N\times N$ matrix, as soon as $j\leq N$ and $k \leq 2j-2$.
When $A =J_\mu$,  the Jacobi matrix associated to a p.m. $\mu$, we get, in terms of the moments :
\ben
\label{cotes}m_k(\mu) = \langle e_1, \left(J_\mu^{[j]}\right)^{k}e_1\rangle \ , \ k=1, \dots, 2j-1.
\een
for every $j$ if $\mu$ as an infinite support, and for $j\leq N$ if $\mu$ is supported by $N$ points.
Notice that this kind of formula leads to Gauss-Jacobi quadratures.
It means that, there exists a sequence of polynomials $f_r$ of $2[N/2]+1$ variables, such that 
\ben
\label{abm}
m_r(\mu) = f_r (b_0, \dots, b_{[r/2]}; a_0, \dots, a_{[r/2]-1})\,,
\een
for any $r$ if $\mu$ as an infinite support, and for $r\leq 2N -1$ if $\mu$ is supported by $N$ points.

Notice that the inverse relations are quite intricated (see for instance Simon \cite{Simonself} Theorem A2). Actually, $a_n$ depends on $m_1, \dots, m_{2n+2}$ and $b_n$ depends on $m_1, \dots, m_{2n+1}$.

\subsection{Tridiagonal representations of $\beta$-ensembles}
We now consider the Jacobi mapping for our random matrix models. The case of the $\beta$-ensembles is directly obtained by
the representation proposed by 
Dumitriu and Edelman (\cite{DuEdel2}). 
\begin{itemize}
\item For the 
normalized G$\beta$E this representation is 
\[H_\beta\sn=\begin{pmatrix} b_0\sn&a_0\sn&0&\dots&0\\
a_0\sn&b_1\sn &a_1\sn&\ddots&\vdots\\
0&\ddots &\ddots&\ddots&0\\
\vdots&\ddots&a_{N-3}\sn&b_{N-2}\sn&a_{N-2}\sn\\
0&\dots&0&a_{N-2}\sn&b_{N-1}\sn
\end{pmatrix}
\]
where the variables $a_0\sn, \dots, a_{N-2}\sn, b_0\sn, \dots, b_{N-1}\sn$ are independent and
\ben\nonumber
b_j\sn &\el& {\mathcal N}(0; (\beta'N)^{-1}) \ \,, \\
\label{loisCCGUE} 
a_j\sn &\el& \sqrt{\gamma\Big(\beta'(N-1-j), (\beta'N)^{-1}\Big)}\,.
\een
It means that $H_\beta\sn$ has the same joint distribution of eigenvalues as for the $G\beta E(N)$. Moreover the weights are independent of the eigenvalues and have the required distribution.

\item For the  L$\beta$E $(N,m(N))$ the representation is $L_\beta\sn = B_\beta\sn \ \big(B_\beta\sn)^T$
\[B_\beta\sn=\begin{pmatrix} d_1\sn&0&0&\dots&0\\
s_1\sn&d_2\sn&0&\ddots&\vdots\\
0&\ddots &\ddots&\ddots&0\\
\vdots&\ddots&s_{m(N)-2}\sn&d_{m(N)-1}\sn&0\\
0&\dots&0&s_{m(N)-1}\sn&d_{m(N)}\sn
\end{pmatrix}
\]
where the variables $d_1\sn, \dots, d_{m(N)}\sn, s_1\sn, \dots, s_{m(N)-1}\sn$ are independent and
\ben
\nonumber
s_j\sn &\el& \sqrt{\gamma\Big(\beta'(m(N)-j), (\beta'N)^{-1}\Big)}\ \,, \\  
\label{loisCCLUE}
d_j\sn &\el& \sqrt{\gamma\Big(\beta'(N+1-j), (\beta'N)^{-1}\Big)}\,.
\een
\item The representation of the $J\beta E(N;a,b)$ has been obtained by Killip and Nenciu (\cite{Killip1}). Actually, they consider a measure $\mu$ on $[-2, 2]$ with finite support as the projection of a symmetric measure $\tilde \mu$ on the unit circle $\mathbb T = \{ z : |z| = 1\}$ by the mapping $z \mapsto z + z^{-1}$. The Jacobi parameters $(a_0, \dots ; b_0, \dots)$ of $\mu$ are in bijection with the Verblunsky coefficients $(\alpha_0, \dots )$ of $\tilde \mu$ by the Geronimus relations (this is also true for measures with infinite support, see \cite{Simon3} section 11)). Notice that choosing a probability distribution to sample Verblunsky coefficients leads to a probability distribution on Jacobi matrices.
\begin{theorem}[Killip-Nenciu, Theorem 2]
\label{kn}
Given $\beta > 0$, let $\alpha_k\sn, 0 \leq k \leq 2N-2$ be independent and distributed as follows:
\[\alpha_{2p}\sn \el\beta_s \Big((N-p-1)\beta' + \a + 1, (N-p-1)\beta' + \b + 1\Big)\ \,,\]
\ben
\label{alphadef}
\alpha_{2p-1}\sn \el\beta_s \Big((N-p-1)\beta'+ \a + \b +2, (N-p)\beta'\Big)\,,
\een
for  $p=0, \dots, N-1$.
Let $\alpha_{2N-1}\sn = \alpha_{-1}\sn = -1$ and define\footnote{these are the Geronimus relations}
\ben
\nonumber
b_k\sn = (1-\alpha_{2k-1}\sn)\alpha_{2k}\sn-(1+\alpha_{2k-1}\sn)\alpha_{2k-2}\sn\\
\label{bdef}
a_k\sn = \sqrt{(1-\alpha_{2k-1}\sn)(1-(\alpha_{2k}\sn)^2)(1+\alpha_{2k+1}\sn)}
\een 
Then the eigenvalues of the tridiagonal matrix $A\sn$ built with these coefficients $a_k\sn$ and $b_k\sn$ 
have a joint density proportional to
\[|\Delta(x_1, \dots, x_N)|^\beta \prod_{j=1}^N (2-x_j)^\a (2 + x_j)^\b\]
and the vector of weights is Dir$(\beta')$ distributed. 
\end{theorem}
We call ${J\beta E}(N, \a, \b)$ ensemble the above distribution on tridiagonal $N\times N$ matrices. 

Since it is often convenient to work on $[0,1]$ instead of $[-2, 2]$, let us introduce the affine mappings : 
\ben\label{affine}x \in [0,1]\stackrel{r}{\mapsto} 4x -2\\
y \in [-2, 2] \stackrel{s}{\mapsto} \frac{y+2}{4}\een
We call $\widehat{J\beta E}(N, \a, \b)$ the image of ${J\beta E}(N, \a, \b)$ by $s$. The 
 preceding result may be rephrased in the following way: 
\begin{corollary}
If $A\sn$ is sampled in the $\widehat{J\beta E}(N, \a, \b)$ ensemble, its eigenvalues 
  have a joint density proportional to
\[|\Delta(x_1, \dots, x_N)|^\beta \prod_{j=1}^N x_j^\b (1- x_j)^\a\,,\]
and the vector of weights is Dir$(\beta')$ distributed. 
\end{corollary}
\end{itemize}
\section{Large Deviations in the $\beta$-Hermite ensemble}
\label{sHerm}
\subsection{Introduction}Recall  that the sequence of ESD
\[
\mu_\u\sn = \frac{1}{N}\sum_{k=1}^N \delta_{\lambda_k}
\]
 satisfies the LDP with speed $\beta'N^2$ and good rate function 
\[I^\u (\mu) = -\Sigma(\mu) + \int_{\mathbb R} \frac{x^2}{2} d\mu(x) + K_H\,,\]
where $K_H$ is a constant (see \cite{BenGui}) and 

\[\Sigma(\mu) = \iint_{\mathbb R^2} \log|x-y| \!\ d\mu(x)d\mu(y).\]
 
The equilibrium measure, unique minimizer of $I^\u$, is the semicircle distribution (denoted hereafter SC, see Section 
\ref{appendixb}). 
In particular, the sequence $(\mu_\u\sn)$ converges weakly in probability to SC.

To manage the large deviations of $(\mu_\w\sn)$, we will first tackle the large deviations of $(a_k\sn , b_k\sn , \ k \geq 0)$. It is important 
to notice already that, in view of (\ref{loisCCGUE}),  as $N \rightarrow \infty$, we have for fixed $k \geq 0$ , $a_k\sn \rightarrow 1$ and $b_k\sn \rightarrow 0$ (in probability). 
The 
 corresponding infinite Jacobi matrix which satisfies
\[b_k = 0 \ , \ a_k = 1 \ , \ \ k\geq 0\,,\]
(often called the free Jacobi matrix, see  Simon (\cite{Simon1} p.13) 
is $J_\mu$ with $\mu = SC$.

In the  large deviations properties of $(\mu_\w\sn)$, 
the extremes eigenvalues will play an important role.
As a matter of fact,  the following function will appear in our rate function.
Let, for $x \geq 2$
\be {\mathcal F}_G (x) = \int_2^{x} \sqrt{t^2 -4}\!\ dt &=& 4 \int_1^{\frac{x}{2}} \sqrt{t^2-1}\!\ dt\\
&=& \frac{x}{2}\sqrt{x^2-4} -2\log\left(\frac{x + \sqrt{x^2 -4}}{2}\right)\,.\ee
Further, for $x<-2$ set ${\mathcal F}_G (x)={\mathcal F}_G (-x)$. The following lemma gives the large deviations properties for the largest eigenvalues in the $G\beta E(N)$ model frame.

\begin{lemma}
\label{LDPmax}
For the $G\beta E(N)$ model the sequence  $(\lambda_{max}\sn)$ satisfies for $x \geq 2$ 
\ben
\lim_N \frac{1}{\beta' N} \log \mathbb P(\lambda_{max}\sn \geq x) = - {\mathcal F}_G(x)\,.
\een 
\end{lemma}
The statement and proof for the GOE are due to \cite{BenAging} Theorem 6.1, the case GUE is in \cite{MR2336602} Prop. 3.1. More generally, for a continuous potential $V$, the result is tackled in  \cite{agz} Theorem 2.6.6. (the potential $V$ in the last theorem is quadratic).

To prepare the statement of our main result, we need another definition.
 \begin{definition}[Simon] 
We say that a p.m. $\mu$ on $\mathbb R$ satisfies the Blumenthal-Weyl condition (B.W.c) if 
\begin{itemize}
\item[i)] $\hbox{Supp}(\mu) = [-2, 2]\cup\{E_j^-\}_{j=1}^{N^-}\cup\{E_j^+\}_{j=1}^{N^+}$ where $N^+$ (resp. $N^-$) is either $0$, 
finite or infinite, 
\[E_1^- < E_2^- < \dots  < -2\ \ \hbox{and}\ \ E_1^+ > E_2^+ > \dots  > 2\] are isolated points of the support.
\item[ii )] If $N^+ = \infty$ (resp. $N^- = \infty$) then $E_j^+$  converges towards $2$ (resp. $E_j^-$ converges towards $-2$).
\end{itemize}
\end{definition}
\subsection{Main result}
Here is our main result. Notice that, of course, SC is the unique minimizer of  the rate function, in accordance with the remark at the beginning of this section.
\begin{theorem}
\label{mainresult}
The sequence  $(\mu_\w\sn)$ satisfies the LDP in ${\mathcal M}^1_{m,d}$ with speed $\beta'N$ and  good rate function
\begin{equation}
\label{projrate}
I^\w(\nu) = 
\begin{cases}  {\mathcal K}(\hbox{SC}\!\ |\!\ \nu) +  \sum_{n=1}^{N^+} {\mathcal F}_G(E_n^+)  +  \sum_{n=1}^{N^-} {\mathcal F}_G(E_n^-)  & \text{if $\nu$  satisfies  B.W.c,}
\\
+\infty &\text{otherwise}\,.
\end{cases}
\end{equation}
\end{theorem}
\noindent
Hence,  the rate function $I(\nu)$ is finite if and only if
\[\nu(dx) = f_a(x) SC(dx) + \nu_s (dx) + \sum_{n=1}^{N^+} \kappa_n \delta_{E_n^+}(dx) + \sum_{n=1}^{N^-} \kappa_n \delta_{E_n^-}(dx) \] 
where  $\nu_s$  is singular (with respect to the Lebesgue measure) and is supported by a subset of $[-2, +2]$ and \[
-\int_{-2}^2 \log f_a (x) SC(dx) < \infty \ , \ 
\sum_{n=1}^{N^+} {\mathcal F}_G(E_n^+)  +  \sum_{n=1}^{N^-} {\mathcal F}_G(E_n^-) < \infty.\]
In this case
\ben I(\nu) = \int_{-2}^2 \log \left(\frac{\sqrt{4-x^2}}{2\pi f_a (x)}\right) \frac{\sqrt{4-x^2}}{2\pi}dx + \sum_{n=1}^{N^+} {\mathcal F}_G(E_n^+)  +  \sum_{n=1}^{N^-} {\mathcal F}_G(E_n^-)\,.\een

\begin{proof}
For $k > 0$, the subset $M(k)$ of ${\mathcal M}_m^1$ of all p.ms supported by $[-k, +k]$ is compact for our topology. 
Indeed, for p.ms in $M(k)$ the moment maps are continuous function ($M(k)$ is tight for the convergence in law).
\noindent
From Lemma \ref{LDPmax} we know that
\[\lim_{k \rightarrow \infty}\lim_{N \rightarrow \infty}\frac{1}{N} \log \mathbb P (\lambda_{max}\sn > k) = - \infty\,.\]
By symmetry, we have also
\[\lim_{k \rightarrow \infty}
\lim_{N \rightarrow \infty}\frac{1}{N}
 \log \mathbb P (\lambda_{min}\sn < -k) = - \infty\,.\]
This implies
\[\lim_{k \rightarrow \infty}\lim_{N \rightarrow \infty}\frac{1}{N} \log \mathbb P ( \mu_\w\sn \notin M(k)) = - \infty\,,\]
hence the sequence $(\mu_\w\sn)_N$ is exponentially tight.
\\
\\
From the inverse contraction principle (see \cite{DZ} Theorem 4.2.4 and Remark a)) it is a consequence of the two following theorems: the first one is a LDP for the sequence of moments and the second one is a {\it magic} formula which allows a powerful identification of the rate function.   
\end{proof}
\noindent
We now give one of the main ingredients of our LDP proof for G$\beta$E(N) ensembles. First define the functions
$g(x) := x- 1 - \log x$ if $x > 0$ and $g(x) = \infty$ otherwise and let
\[G(x) := \begin{cases}
       g(x^2)  & \text{if $x > 0$} \\
      \infty & \text{otherwise}\,.
     \end{cases}\]
\begin{theorem}
\label{LDPmom}
The sequence  $\big({\mathbf m}(\mu_\w\sn)\big)$  satisfies in $\mathbb R^{\mathbb N}$ the LDP with speed $\beta' N$ and good rate function $I$ defined as follows. This function is finite if and only if there exist $(b_0, \dots ; a_0, \dots)\in \mathbb R^{\mathbb N}\times (0, \infty)^{\mathbb N}$ 
satisfying
\ben
\label{cpct}
\sum_{j=0}^\infty b_j^2 < \infty \ , \ \  \sum_{j=0}^\infty (a_j - 1)^2 < \infty
\een
such that 
$m_r = \langle e_1, A^r e_1\rangle$ for every $r \geq 1$ with $A$ infinite tridiagonal matrix built with  $(b_0, \dots ; a_0, \dots)$. 
In that case \[I(m_1, \dots) = \frac{1}{2}\sum_{j=0}^\infty b_j^2 + \sum_{j=0}^\infty G(a_j)< \infty\,.\]
\end{theorem}

\begin{theorem}[Killip-Simon \cite{Sum1}, \cite{Simon2} Theorem 13.8.6]
Let $J$ be a Jacobi matrix built with $(a_0, \dots; b_0, \dots)\in (0, \infty)^{\mathbb N} \times {\mathbb R}^{\mathbb N}$ satisfying $\sup a_n + \sup |b_n| < \infty$.  Let $\mu$ be the associated measure obtained by Favard's theorem.
Then
\ben
\sum_k\left[ b_k^2 + (a_k - 1)^2\right] < \infty
\een
if, and only if, the p.m. $\mu$ satisfies  B.W.c. and the two following conditions:
\ben
\sum_{j=1}^{N_+} (E_j^+ -2)^{3/2}+ \sum_{j=1}^{N_-} (-2 - E_j^-)^{3/2}< \infty\\
\int_{-2}^2 \log (f_a (x))\sqrt{4-x^2}\!\ dx > -\infty\,.
\een
In that case
\ben
I^\w(\mu) = \sum_n \left[\frac{1}{2} b_n^2 + G(a_n)\right]
\een
where both sides may be (simultaneously) infinite.
\end{theorem}
The proof of Theorem \ref{LDPmom} will use the following result.

\begin{lemma}
\label{LDPab}
For fixed $k$, 
 $\big(b_0\sn, \dots, b_k\sn; a_0\sn, \dots, a_{k-1}\sn\big)_{N \geq k}$ satisfies
in $\mathbb R^{2k-1}$ a LDP with speed $\beta'N$ and good rate function
\ben I_k(b_0, \dots, b_k; a_0, \dots, a_{k-1}) = \frac{1}{2}\sum_{j=0}^k b_j^2 + \sum_{j=0}^{k-1} G(a_j)\,.\een
\end{lemma}
\begin{proof}
 It is an immediate consequence of independence and  the LDP for Gaussian and Gamma r.vs. recalled in the following lemma.
\end{proof}
\begin{lemma}
\label{PGDclass}
\ 
\begin{enumerate}
\item The sequence of distributions ${\mathcal N}(0; n^{-1})$ satisfies the LDP with speed $n$ and good rate function $x \mapsto x^2/2$.
\item
For $\alpha > 0$ and $c$ fixed, the sequence of distributions $\gamma\big((n- c), (\alpha n)^{-1}\big)$ satisfies the LDP with speed $n$ and good rate function $x \mapsto g(\alpha x)$.
\item
\label{PGDbetas}
For $u,v >0$ and $\delta, \delta'$ fixed, the sequence of distributions 
$\beta_s(un +\delta, vn + \delta')$ satisfies the LDP with speed $n$ and good rate function:
\ben
\label{hq}
 h(q)   =\begin{cases} q(u-v) -u \log(1+q) -v \log(1-q) \;\; ; \;\; q \in (-1, 1)\\
\infty\;\;\ ; \;\; \mbox{otherwise} \,.
\end{cases}
\een 
\end{enumerate}
\end{lemma}
\begin{proof}
The points 1 and 2 are well known. For point 3, we use the representation
\[\beta_s(un +\delta, vn + \delta') \el \frac{\gamma(un +\delta) - \gamma(vn +\delta')}{\gamma(un +\delta) + \gamma(vn +\delta')}\]
hence by contraction the rate function is
\[h(q) = \inf\{ ug(x/u)+vg(y/v); \frac{x-y}{x+y} = q\}\,, 
\]
which yields easily (\ref{hq}).
\end{proof}

\begin{proof}{{\bf of Theorem \ref{LDPmom}.}}
Fix $\ell>1$. By Lemma \ref{LDPab} and the contraction principle, the sequence
$\big(m_1(\mu_\w\sn) , \dots, m_{2\ell-1}(\mu_\w\sn)\big)$ satisfies the LDP in $\mathbb{R}^{2\ell-1}$ with speed $\beta'N$ and rate function $\wt I_{2\ell-1}$ defined as follows. Notice that there is at most only one tridiagonal matrix $A_\ell$ built from $(b_0, \dots, b_{\ell-1} ; a_0, \dots,a_{\ell-2})$
as in (\ref{favardfini}) such that
\begin{equation}
m_r = \langle e_1, A_\ell^r e_1\rangle , \ \ r=1, \dots, 2\ell-1.
\label{chavez}
\end{equation}
Hence, 
if $(m_1, \dots, m_{2\ell-1})$ satisfies (\ref{chavez}), then 
\ben
\wt I_{2\ell -1}( m_1, \dots, m_{2\ell-1})&=& I_{\ell-1}(b_0, \dots, b_{\ell-1} ; a_0, \dots,a_{\ell-2})
\een
Otherwise, $\wt I_{2\ell-1}(m_1, \dots, m_{2\ell -1})$ is infinite.
\noindent
We do not consider the even case since there is no injectivity in that case.

\noindent
We now apply the Dawson-G\"artner theorem. 
 Let us endow $\mathbb{R}[[X]]$ with the topology of pointwise convergence of coefficients. It can be viewed  as the projective limit
 \[\mathbb{R}[X] = \lim_{\longleftarrow}\mathbb{R}_k[X]\] where $\mathbb{R}_k[X]$
is the set of polynomials of degree equal or less than $k$.
\\
\noindent
The rate function is
\ben
I(m_1, \dots) = \sup \{\wt I_{2k+1}(m_1, \dots, m_{2k+1}) : k \geq 0\}\,.
\een
It is clear that \be\sup \{\wt I_{2k+1}(m_1, \dots, m_{2k+1}) : k \geq 0\} &=& \sup_k \{\frac{1}{2}\sum_{j=0}^k b_j^2 +
\sum_{j=0}^{k-1} G(a_j)\}\\
&=& \frac{1}{2}\sum_{j=0}^\infty b_j^2 +
\sum_{j=0}^\infty G(a_j)\leq \infty\,.\ee
\end{proof}
\subsection{Failure of the direct method}

Mimicking the unitary case (\cite{FABALOZ}), it is tempting to define the random measure 
\[\wt\mu_\w\sn = \sum_{k=1}^N Y_k \delta_{\lambda_k}\]
with the $Y_k$ independent and $\gamma(\beta')$ distributed so that
\[\mu_\w\sn = \frac{\wt\mu_\w\sn}{\wt\mu_\w\sn(1)}\]
The problem is that the general method of Najim \cite{Najim1} cannot be applied. Indeed, the main assumption on the range of the  eigenvalues is violated. As a matter of fact, not all the eigenvalues belong to the support of the semicircle law. Outliers give a contribution.
So that,  the conclusion given by this approach is not true. The rate function candidate only contains the Kullback part of the LDP but loose the {\it outer} part.
\section{Large Deviations in the $\beta$-Laguerre ensemble}
\label{sLag}
In the Laguerre case, in the usual asymptotics $N \rightarrow \infty$, $m(N)/N \rightarrow \tau < 1$, we observe  similar phenomena. 
Recall that the sequence of ESD
\[\mu_\u\sn = \frac{1}{m(N)}\sum_{k=1}^{m(N)} \delta_{\lambda_k}\]
satisfies the LDP with speed $\beta'N^2$ and good rate function
\[I^\u (\mu) = -\tau^2\Sigma(\mu) +\tau\int_0^\infty \Big(\frac{x}{2} - (1-\tau)\log x\Big)\!\ d\mu(x) + K_L\,,\] 
where $K_L$ is a constant (\cite{hiai1}).
%
 The equilibrium measure, unique minimizer of $I^\u$  is the Marchenko-Pastur distribution of parameter $\tau$ (denoted hereafter by MP, see Appendix). In particular, the sequence $(\mu_\u\sn)$ converges weakly in probability to MP.
 
 To manage the large deviations of $(\mu_\w\sn)$, we will first tackle the large deviations of $(s_k^{N,m} , d_k^{N,m} , \ k \geq 0)$. 
Recall that the elements of the tridiagonal matrix $L_\beta\sn$ are
\ben\nonumber
 b_0\sn = (d_1\sn)^2 \ ,\ b_k\sn  &=& (s_k\sn)^2 + (d_{k+1}\sn)^2 \ \ (1\leq k \leq N-1) \\
 \label{dsdonneabN}
a_k\sn &=&  s_{k+1}\sn d_{k+1}\sn \ \ \ \ \ (0 \leq k \leq N-2)\,.\een
We can see already that, in view of (\ref{loisCCLUE}), we have for fixed $k \geq 1$ and $N \rightarrow \infty$, 
$\lim d_k\sn = 1$ and $\lim s_k\sn = \sqrt \tau$
 (in probability). From (\ref{dsdonneabN}), this yields 
$\lim b_0\sn = 1$  
and for fixed $k \geq 1$, $\lim b_k\sn = 1+\tau$ , $\lim a_{k-1}\sn = \sqrt \tau$ (in probability). 
 The 
 corresponding infinite  Jacobi matrix which satisfies
\[b_0 = 1\ , \ b_k = 1 +\tau \ , \ (k\geq 1)\ \ ; \ \ a_k = \sqrt\tau \ ; \ (k \geq 0)\,.\]
is $J_\mu$ with $\mu = MP$.

Let ${\mathcal F}_L$ defined by
\[{\mathcal F}_L(x) = \begin{cases} \displaystyle\int_{b(\tau)}^x \frac{\sqrt{(t -  a(\tau))(t - b(\tau))}}
{t\tau}\!\ dt \;\;\;\;x \geq b(\tau)\,,\\ \displaystyle
\int_x^{a(\tau)} \frac{\sqrt{(a(\tau)-t)(b(\tau) -t)}}
{t\tau}\!\ dt        \;\;\;\; 0 < x \leq a(\tau)\,. 
\end{cases}\]
\begin{lemma}
\label{lemsale}
For the L$\beta$E$(N, \tau N)$ model, 
\begin{enumerate}
\item
the sequence  $(\lambda_{max}\sn)$ satisfies for $x\geq b(\tau)$
\ben
\lim_N \frac{1}{\beta' N} \log \mathbb P(\lambda_{max}\sn \geq x) = - {\mathcal F}_L(x)
\,.
\een
\item 
the sequence  $(\lambda_{min}\sn)$ satisfies for $0 < x \leq a(\tau)$
\ben
\lim_N \frac{1}{\beta' N} \log \mathbb P(\lambda_{min}\sn \leq x) = - {\mathcal F}_L(x)
\,.\een
\end{enumerate}
\end{lemma}
\begin{remark}
\label{r52}
As mentioned before,  a LDP for a general continuous potential is proved in \cite{agz} Theorem 2.6.6. The knowledge of the Cauchy-Stieltjes transform of $MP$ allows to recover the formula given in \cite{Feral} p. 47. Here, the potential is
$$V(x)=\tau\frac{x}{2} -\tau (1-\tau)\log x.$$
\end{remark}

For a general double sequence of positive numbers $(d_k)_{k \geq 1}$ and $(s_k)_{k\geq 1}$ we set
 $d\circ s = (d_1, \dots ; s_1, \dots )$. We deduce the elements 
 \ben
\nonumber
 b_0 = d_1^2 \ ,\ b_k  = s_k^2 + d_{k+1}^2 \ \ (k \geq 1) \\
 \label{dsdonneab}
a_k  = s_{k+1}d_{k+1} \ \ (k \geq 0)\,.\een
Conversely, if $(a_0, \dots ; b_0, \dots)$ is given in $(0,\infty)^{\mathbb N \times\mathbb N}$ such that the tridiagonal matrix is positive, there exists a unique $d\circ s$ satisfying (\ref{dsdonneab}).
Here is a direct consequence of Lemma \ref{PGDclass}.

\begin{theorem}
\label{5.2}
Under the $L\beta E(N, \tau N)$ model, the sequence  $(\mu_\w\sn)$ satisfies in ${\mathcal M}_{m,d}^1((0, \infty))$ a LDP with speed $\beta'N$ and good rate function $I^\w$ defined as follows. This function is finite at $\nu$ if and only if there exists $d\circ s \in [0, \infty)^{\mathbb N}\times [0, \infty)^{\mathbb N}$ (necessarily unique) satisfying
\[\sum_k G(d_k) + \tau \sum_k G(s_k/\sqrt{\tau}) < \infty\,,\]
such that $m_r(\nu) =\langle  e_1, A^r e_1\rangle$ for every $r \geq 1$ with $A$ infinite tridiagonal matrix built with $(b_0, \dots ; a_0, \dots)$ satisfying (\ref{dsdonneab}). In that case
\ben\label{sum0}I^\w(\nu) = \sum_k G(d_k) + \tau \sum_k G(s_k/\sqrt{\tau})\,.\een
\end{theorem}

\begin{remark}
\
\begin{itemize}
\item It is clear from (\ref{sum0}) that the unique minimizer of $I^\w$ corresponds to $d_k \equiv 1$ and $s_k \equiv \sqrt \tau$  which corresponds to MP.
\item 
When $\tau=1$, we can write:
\be I^\w(\nu) &=& \sum_{k\geq 1} \big[d_k^2 - 1 - \log d_k^2 + s_k^2 - 1 - \log s_k^2\big] 
\\ &=& d_1^2 -1 + \sum_{k\geq 1} \big[d_{k+1}^2 +s_k^2 - 2\big] - 2\sum_{k\geq 1}  \log (d_ks_k) \\
 &=& b_0 - 1 +\sum_{k\geq 1} (b_k - 2) - 2 \sum_{k\geq 0} \log a_k\,.
\ee
This expression of $I^\w$ in terms of the Jacobi coefficients makes plausible the existence of a convenient sum rule and we propose the following conjecture : 
\end{itemize}
\end{remark}
\begin{conj}
The rate function is
\[I^\w(\nu) = {\mathcal K}( MP\!\ |\!\ \nu) + \sum_j {\mathcal F}_L(E_j^\pm)\,.\] 
\end{conj}
\begin{proof}{\bf of Theorem \ref{5.2}} 
For $k$ fixed, $(d_k\sn)$  (resp. $(s_k\sn$))  satisfies a LDP with good rate function $G(d_k)$ (resp. $\tau G(s_k/\sqrt{\tau})$) hence, by independence, the rate function is the sum (\ref{sum0}).
 \end{proof}
\section{Large Deviations in the $\beta$-Jacobi ensemble}
\label{sJac}
Let us consider the $\widehat{J\beta E}(N, \a(N), \b(N))$ ensemble. The usual asymptotics  is $N \rightarrow \infty$, $\b(N)/N \rightarrow \beta'\kappa_1$, $\a(N)/N \rightarrow \beta'\kappa_2$. The sequence of ESD
\[\mu_N = \frac{1}{N}\sum_{k=1}^N \delta_{\lambda_k}\]
satisfies the LDP with speed $\beta'N$ and good rate function :
\ben\label{entrjac} I^\u (\mu) = -\Sigma(\mu) -\int_0^1 (\kappa_1 \log x + \kappa_2 \log(1-x))\!\ d\mu(x) + K_J\,,\een
where $K_J$ is a constant (see \cite{hiai2}). 
The equilibrium measure, unique minimizer of  $I^\u$ is the Kesten-MacKay distribution (denoted hereafter KMK) of parameter $(u_- , u+)$, 
where 
\[u_-, u_+ = u_\pm \left(\frac{1+\kappa_1}{2+\kappa_1+\kappa_2}, \frac{1+\kappa_1+\kappa_2}{2+\kappa_1+\kappa_2}\right)\]
(see Section \ref{appendixb}). 
In particular, the sequence $(\mu_\u\sn)$ converges weakly in probability to KMK.

To manage the large deviations of $(\mu_\w\sn)$, we will first tackle the large deviations of $(\alpha_k\sn , k \geq 0)$. It is important to notice already that, in view of (\ref{alphadef}), we have for fixed $p \geq 0$, \[\lim_N \alpha_{2p}\sn = \frac{\kappa_2 - \kappa_1}{2 + \kappa_1 + \kappa_2}\ , \ \lim_N \alpha_{2p+1}\sn = -\frac{\kappa_1 + \kappa_2}{2 + \kappa_1 + \kappa_2}\,.\]

The symmetric measure admitting these limiting Verblunsky coefficients is well understood by its Cauchy-Stieltjes transform since the work of Geronimus (\cite{Gero1}, see also the books of Simon). We do not give details here to shorten the paper. After projection, we obtain the KMK distribution.

Let ${\mathcal F}_J$ defined by
\[{\mathcal F}_J(x) = \begin{cases} \displaystyle\int_{u_+}^x \frac{\sqrt{(t -  u_+)(t - u_-)}}
{t(1-t)}\!\ dt \;\;\;\; u_+ \leq  x < 1\,,\\ \displaystyle
\int_x^{u_-} \frac{\sqrt{(u_--t)(u_+ -t)}}
{t(1-t)}\!\ dt       \;\;\;\; 0 < x \leq u_-\,. 
\end{cases}\]
This following lemma is a kin of  Lemmas \ref{LDPmax} and  \ref{lemsale}. Here, the potential is 
\[V(x) = - \kappa_1 \log x - \kappa_2 \log (1-x)\,.\]
\begin{lemma}
For the $\widehat{J\beta E}(N, \a(N), \b(N))$ model with the above notations, if $\kappa_1, \kappa_2 > 0$,
\begin{enumerate}
\item
the sequence  $(\lambda_{max}\sn)$ satisfies for $x\in (u_+, 1)$
\ben
\lim_N \frac{1}{\beta' N} \log \mathbb P(\lambda_{max}\sn \geq x) = - {\mathcal F}_J(x)
\,.\een
\item 
the sequence  $(\lambda_{min}\sn)$ satisfies for $x \in (0, u_-)$
\ben
\lim_N \frac{1}{\beta' N} \log \mathbb P(\lambda_{min}\sn \leq x) = - {\mathcal F}_J(x)
\,.\een
\end{enumerate}
\end{lemma}


\begin{theorem}
\
\begin{enumerate}
\item (Gamboa-Rouault \cite{FABALOZ})
 Under the $\widehat{J\beta E}(N, \a,\b)$ model, the sequence  $(\mu_\w\sn)$ satisfies in ${\mathcal M}^1([0, 1])$ endowed with the weak topology the LDP with speed $N$ and  good rate function \[I(\nu) =  {\mathcal K}(ARCSINE\!\ |\!\ \nu)\,.\]
\item
Under the $\widehat{J\beta E}(N, \kappa_1 N, \kappa_2 N)$ model, the sequence  $(\mu_\w\sn)$ satisfies in ${\mathcal M}^1([0, 1])$ endowed with the weak topology the LDP with speed $N$ and with a good rate function $I^\w$ defined as follows. This function is finite at $\nu$ if and only if there exists $\vec{\alpha}\in (-1, 1)^{\mathbb N}$ (necessarily unique) such that 
\be {\mathcal I}(\vec{\alpha})&:=&  (\kappa_1 - \kappa_2) \sum_0^\infty \alpha_{2k} + (\kappa_1 + \kappa_2) \sum_0^\infty \alpha_{2k+1}\\
&-& (1+\kappa_1) \sum_0^\infty \log(1+\alpha_{2k}) - (1+\kappa_2) \sum_0^\infty \log(1-\alpha_{2k})\\ &-& (1+\kappa_1 +\kappa_2) \sum_0^\infty \log(1+\alpha_{2k+1}) -  \sum_0^\infty \log(1-\alpha_{2k+1})\ee
is finite. Here $\vec{\alpha}$ is related to $\nu$ through Geronimus relation (see \ref{bdef}).
In that case 
\[I^\w(\nu) = {\mathcal I}(\vec{\alpha})\,.\]
\end{enumerate}
\end{theorem} 
\begin{proof}
We apply Lemma \ref{PGDclass} (3), with $n =\beta'N$, and for an even index we have $u = 1+\kappa_1 , v = 1 + \kappa_2$ and with odd index 
$u= 1+\kappa_1 + \kappa_2 , v = 1$
\be I_{\alpha_{2k}}(x) &=&  x(\kappa_1 - \kappa_2) - (1+\kappa_1)\log (1+x) - (1+ \kappa_2)\log(1-x)\\
I_{\alpha_{2k+1}}(x) &=& x(\kappa_1 + \kappa_2) - (1+\kappa_1+\kappa_2)\log (1+x) - \log(1-x)
\ee
Then it is enough to add up.

In the particular case of $\a$ and $\b$ fixed, we have $\kappa_1 = \kappa_2 = 0$ and 
\[I(\vec{\alpha}) = - \sum_0^\infty \log (1 - \alpha_k^2)\,.\]
But the Szeg\"o formula (\cite{Simon1}) says that it is exactly the reversed Kullback with respect to the ARCSINE distribution.
\end{proof}

In the general case, there is up to our knowledge, no known sum rule. Besides it is very intricate to express the above sums in terms of the tridiagonal coefficients. Nevertheless it is tempting to propose the conjecture.

\begin{conj}
Under the $\widehat{J\beta E}(N, \kappa_1 N, \kappa_2 N)$ model, the rate function is given by
\[I(\nu) = {\mathcal K}(KMK\!\ |\!\ \nu) + \sum_j {\mathcal F}_J(E_j^\pm)\]
\end{conj}
\subsection{Some distributions}
\label{appendixb}
\subsubsection{Gamma distribution}
For $a,b>0$, the $\gamma(a,b)$ distribution is supported by $[0, \infty)$ with density
\[\frac{e^{-x/b} x^{a-1}}{b^a \Gamma(a)}\]
Its mean is $ab$.
\subsubsection{Beta distribution}
For $a,b>0$, the beta symmetric distribution of parameter $(a,b)$, denoted by $\beta_s(a,b)$, is supported by $(-1, 1]$ and has  density
\[2^{1-a-b}\frac{\Gamma(a+b)}{\Gamma(a)\Gamma(b)}(1-x)^{a-1} (1+x)^{b-1}\]
Its mean is $\frac{b-a}{b+a}$.
\subsubsection{Dirichlet distribution}
For $k \geq  1$, we set
\be {\mathcal S}_{k} := \{(x_1, \cdots, x_{k}) : x_i > 0 , (i = 1, \cdots , k) ,  \ x_1 + \cdots + x_{k}  = 1 \} \\
{\mathcal S}_k^< := \{(x_1, \cdots, x_k) : x_i > 0 , (i = 1, \cdots , k) ,  \ x_1 + \cdots + x_k  < 1 \}\,.
\ee
Obviously, the mapping $(x_1, \cdots , x_{k+1}) \mapsto (x_1, \cdots, x_k)$ is a bijection from the simplex ${\mathcal S}_{k+1}$ onto  ${\mathcal S}_k^<$.

For $a_j>0,\;j=1,\ldots, k+1$,
the Dirichlet distribution $\hbox{Dir}(a_1, \cdots, a_{k+1})$ on ${\mathcal S}_{k+1}$ 
has the density
\ben
\label{defdir}
\frac{\Gamma(a_1 + \cdots + a_{k+1})}{\Gamma(a_1) \cdots \Gamma(a_{k+1})}\!\ x_1^{a_1 - 1} \cdots x_{k+1}^{a_{k+1} -1}
\een
with respect to the Lebesgue measure on ${\mathcal S}_{k+1}$. When $a_1= \cdots = a_{k+1} = a>0$, we will denote the Dirichlet distribution by
$\hbox{Dir}_k (a)$. If $a=1$ we recover  the uniform distribution on ${\mathcal S}_k^<$.
\subsubsection{Semicircle distribution}
The semicircle distribution denoted by $SC$ is supported by $[-2, 2]$ with density
\[\frac{\sqrt{4-x^2}}{2\pi}\,.\]
Its Cauchy-Stieltjes transform\footnote{Throughout, all branches of the square roots are taken in accordance with the definition of Cauchy transform}
 is
\ben
m(z) = \int \frac{d\mu(x)}{x-z} = \frac{-z + \sqrt{z^2 -4}}{2} 
\een

When $0 < \tau \leq 1$, the Marchenko-Pastur distribution, denoted by $MP$ is supported by $(a(\tau), b(\tau))$ where $a(\tau) = (1 -
\sqrt{\tau})^2$  , $b(\tau) = (1 + \sqrt{\tau})^2$ 
with  density
\ben \label{defpi}
\frac{\sqrt{(x -  a(\tau))( b(\tau) - x)}}
{2\pi \tau x}\,.\een
Its Cauchy-Stieltjes transform is
\ben
m(z) 
 = \frac{-z -1+\tau + \sqrt{(z-1-\tau)^2 - 4\tau}}{2\tau z}\,.
\een
 \subsubsection{Kesten-McKay distribution}
The Kesten-McKay distribution  is  supported by $(u_- , u_+)$ with $0\leq u_- <u_+\leq 1$ and its  density is
\ben
C_{u_-,u_+}\frac{\sqrt{(x-u_-)(u_+ -x)}}{2\pi x (1-x)}
\een
where 
$$C_{u_-,u_+}^{-1}:=\frac{1}{2}\left[1-\sqrt{u_-u_+}-\sqrt{(1-u_-)(1-u_+)}\right].$$
To express its Cauchy-Stieltjes transform, let us give some notation.
For $(b,c) \in (0,1)\times(0,1)$ we put \ben
\label{defsig} \sigma_\pm (b,c) = \frac{1}{2}\left[1 + \sqrt{bc} \pm
\sqrt{(1- b)(1- c)}\right]\,, \een and for $(x,y) \in
(0,1)\times(0,1)$ \ben \nonumber
u_\pm(x,y) &=& (1 -x-y + 2xy) \pm 2 \sqrt{x(1-x)y(1-y)}\\
\label{defla} &=& \left(\sqrt{(1-x)(1-y)} \pm \sqrt{xy}\right)^2\,.
\een The mappings $\sigma_\pm$ and $u_\pm$ are inverse in the
following sense : \ben\label{lets}\{(b,c) : 0 < b < c < 1\}
\xrightleftharpoons[(u_-, u_+)]{(\sigma_- , \sigma_+)} \{(x,y) : 0<
x < y<1 \ \hbox{and} \ x+y > 1\}\een
The Cauchy-Stieltjes
 is then
  (see for instance \cite{Dem} p.129, or \cite{capitaine} p.425) 
\ben
\label{Dem}
m(z) 
= \frac{(1-\sigma_+ - \sigma_-)}{2(1 - \sigma_+)z} + \frac{(\sigma_+ - \sigma_-)}{2(1-\sigma_+)(1-z)} +\frac{\sqrt{(z-a_-)(z-a_+)}}{2z(1-z)}\,.
\een
ARCSINE corresponds to  $u_-=0$ and $u_+=1$.
\\
\\
{\bf Acknowledgment}
Many thanks are due to Professor Holger Dette for helpful discussions.
\bibliographystyle{plain}
\bibliography{toulz}
\end{document}